\documentclass[preprint,12pt]{elsarticle}

\usepackage{amssymb}

\usepackage{amsmath}
\journal{Fractional Calculus and Applied Analysis }

\begin{document}

\begin{frontmatter}

\title{Linear Canonical Jacobi-Dunkl Transform: Theory and 
	Applications  }

\author{Rong-Qian, Linghu$^a$, Bing-Zhao, Li$^{a,b,*}$ 
	\affiliation{organization={School of Mathematics and Statistics},
		addressline={Beijing Institute of Technology}, 
		city={Beijing},
		postcode={102488}, 
		country={P.R. China}}
	\affiliation{organization={Beijing Key Laboratory on MCAACI},
		addressline={Beijing Institute of Technology}, 
		city={Beijing},
		postcode={102488}, 
	    country={P.R. China
	    \\ *$Corresponding author: li_{-}$$bingzhao@bit.edu.cn$}}	    	
	}
\begin{abstract}
This paper aims to develop an innovative method for harmonic analysis by introducing the linear canonical Jacobi-Dunkl transform (LCJDT), which integrates both the Jacobi-Dunkl transform (JDT) and the linear canonical transform (LCT). Firstly, the kernel function of the LCJDT is derived, and its fundamental properties are examined. Subsequently, the LCJDT is established, along with an investigation of its essential properties, including the inversion formula, Parseval's theorem, differentiation, the convolution theorem, and the uncertainty principle. Finally, the potential application of the LCJDT in solving the heat equation is explored.

\end{abstract}

\end{frontmatter}

\section{ Introduction}
Classical integral transforms, including the Fourier transform (FT) \cite{r1}, Laplace transform (LT), and wavelet transform (WT) \cite{r9}, are foundational tools in contemporary mathematics, physics, and engineering disciplines fileds. These transforms enable the conversion of problems from the time domain to the transform domain, thereby facilitating information analysis and problem solving. However, when non-stationary and nonlinear signals are encountered, or when mathematical problems involve complex structures and weighted functions, classical integral transforms often prove inadequate, leading to significant difficulties and challenges \cite{r40}.

In response to these challenges, researchers have introduced new transform methods, such as the Jacobi-Dunkl transform (JDT) \cite{r15} and the fractional Fourier transform (FrFT) \cite{r1}, to more effectively address complex information and mathematical problems. Examples include the fractional Jacobi transform \cite{r17}, fractional Fourier-Jacobi transform \cite{r12}, and fractional Jacobi-Dunkl transform (FRJDT) \cite{r14, r16}, all of which have been proposed to tackle these intricate challenges.

Although the aforementioned researchers have combined fractional operators with the JDT to provide an effective method for handling non-stationary and nonlinear signals, as well as addressing mathematical problems involving complex structures and weighted functions, limitations persist in certain cases \cite{ r12}. Specifically, these methods may not be well-suited for dealing with asymmetry, complex boundary conditions, and partial differential equations (PDE) involving weighted functions \cite{r14, r16}. Moreover, with the increasing demand for advanced information processing, particularly in the directional feature representation of image signals, there is a growing need to capture local signal features at different scales, directions, and positions \cite{r13}. This demand further restricts the applicability of the aforementioned methods. As a result, the search for novel and more effective mathematical approaches to address these challenges has become a focal point of current research \cite{r12, r14, r16, r13}.

The linear canonical transform (LCT) is a linear integral transform with three degrees of freedom \cite{r18}, capable of encompassing classical integral transforms such as the FT \cite{r1}, FrFT \cite{r2}, Fresnel transform (FRT) \cite{r3}, and Lorentz transform (LRT) \cite{r4} under specific conditions. The LCT is widely regarded as a powerful tool in fields such as mathematics, optics, and information analysis \cite{r5, r6, r7, r8}. In this paper, we combine the advantages of the JDT and LCT to propose a novel transform method, the linear canonical Jacobi-Dunkl transform (LCJDT). The LCJDT introduces five degrees of freedom, providing enhanced flexibility and adaptability. It offers strong directional analysis capabilities for capturing local information from multiple angles, thereby improving efficiency and performance in areas such as image encryption and target detection. Furthermore, it effectively addresses PDEs involving asymmetry, complex boundary conditions, and weighted functions \cite{r12, r14, r16, r13, r5, r6, r7}. As a new mathematical tool, the LCJDT demonstrates broad potential for application and provides robust support for solving complex mathematical problems.

The remainder of this paper is structured as follows: In Section 2, some basic preliminaries are provided. In Section 3, the kernel function of the LCJDT is derived, and several basic properties of this kernel function are investigated. In Section 4, the theory of the LCJDT is studied, including its definition, inversion formula, Parseval's identity, basic properties. In Section 5, the uncertainty principles of the LCJDT is studied. In Section 6, the potential application of the LCJDT in solving PDE are explored. In Section 7, the paper is concluded.

\section{ Preliminaries}

In this section, a brief overview of the JDT and LCT is provided, along with the key properties that are essential for developing the main results presented in this paper.

\subsection{\textbf{Jacobi-Dunkl transform (JDT)}}

In this paper, the operator of the JDT is denoted by $\Lambda_{\alpha,\beta}$. It is a partial differential operator that extends the classical Jacobi differential operator by incorporating reflection symmetry \cite{r13}.

\textbf{Definition 2.1:} Let $f(x)\in C^1(\mathbb{R})$, and the parameter $\alpha>-\frac12$, $\beta\geq-\frac12$, $\alpha>\beta$, the $\Lambda_{\alpha,\beta}$ is defined by

\begin{equation}
	\label{a2}
	\Lambda_{\alpha,\beta}f(x)=\begin{cases}\frac{d}{dx}f(x)+\frac{A'_{\alpha,\beta}(x)}{A_{\alpha,\beta}(x)}\frac{f(x)-f(-x)}{2},\text{if}~x~\neq0, \\(2\alpha+2)f'(0), \text{if} ~x~=0,\end{cases}
\end{equation}
where  $A_{\alpha,\beta}(x)=2^{2\rho}\sinh^{2\alpha+1}(|x|)\cosh^{2\beta+1}(x),  \rho=\alpha+\beta+1$.
\\

 \textbf{Proposition 2.1:} Let $f(x)\in C^1(\mathbb{R})$, the Jacobi-Dunkl kernel function  is denoted by $\psi_{\lambda}^{\alpha,\beta}$, which is a unique solution for the following PDE
\begin{equation}
		\label{a3}
\begin{cases} \Lambda_{\alpha,\beta}f(x)=i\lambda f(x),\\ f(0) = 1, \end{cases}
\end{equation}
by solving (\ref{a3}), we obtain 
\begin{equation}
		\label{a4}	
			\left. \psi_{\lambda}^{\alpha,\beta}(x)=\left\{\begin{matrix}\phi_{\mu}^{(\alpha,\beta)}(x)-(i\lambda^{-1} )\phi_{\mu}'^{(\alpha,\beta)}(x), ~if ~ \lambda\neq 0, \\1, ~~if~ \lambda=0, \end{matrix}\right.\right.
\end{equation}
and
\begin{equation}
		\label{a6}
	\phi_\mu^{(\alpha,\beta)}(x)= {}_2F_1\left(\frac{\rho+i\mu}{2},\frac{\rho-i\mu}{2},\alpha+1,-(\sinh x)^2\right),
\end{equation}
where $ \mu^2=\lambda^2-\rho^2$, the $_2F_1$ is the Gauss hypergeometric function \cite{r19,r21}. 
\\

 \textbf{Proposition 2.2:} Throughout the rest of this paper, the following notations will be used, the space of measurable functions $f(x)$ on $\mathbb{R}$, then
\\
 
(i) $L_{\alpha,\beta}^p(\mathbb{R})$, when $p\in[1,\infty)$,   $\begin{aligned}\|f\|_{\alpha,\beta,p}&=\left(\int_{\mathbb{R}}|f(x)|^p\:A_{\alpha,\beta}dx\right)^{1/p}<\infty,\end{aligned}$

(ii) $L_\sigma^p(\mathbb{R})$, when $p\in[1,\infty)$, $\|f\|_{p,\sigma}=\left(\int_\mathbb{R}|f(x)|^pd\sigma(x)\right)^{\frac1p}<+\infty,$
where $d\sigma(x)$ is the measure, the specific mathematical definition can be found in the references \cite{r15},

(iii) $\mathcal{D}(\mathbb{R})$,  $f(x) \in$ $\mathcal{D}(\mathbb{R})$ are smooth (infinitely differentiable) everywhere on $\mathbb{R}$ and are nonzero only within a limited region of the real line. 
\\

\textbf{Definition 2.2:} 
The JDT of a function $f\in L_{\alpha,\beta}^1(\mathbb{R})$ is defined by
\begin{equation}\mathcal{F}_{\alpha,\beta}(f)(\lambda)=\int_\mathbb{R}f(x)\psi_\lambda^{\alpha,\beta}(x)A_{\alpha,\beta}(x)dx.\end{equation}

The JDT combines classical integral transforms with the Dunkl operator, extending the application scope of classical Fourier, Laplace, and Jacobi transforms. It exhibits distinct advantages, particularly when dealing with problems involving reflection symmetry and special structures \cite{r12,r14,r16,r13}. For detailed derivation, definition, and properties of the JDT, please refer to the reference \cite{r15}.

\subsection{\textbf{ Linear canonical transform (LCT)}}

The LCT is a highly versatile and robust mathematical framework, offering a comprehensive approach to the analysis and processing of signals and functions. Its broad applicability spans a diverse range of disciplines, including signal processing, optics, quantum mechanics, and communications, where it serves as an indispensable tool \cite{r18,r30}. The LCT's distinguishing feature lies in its capacity to unify and generalize multiple signal transformations, facilitating the execution of a wide array of operations such as time-frequency analysis, filtering, and system modeling \cite{r6,r7,r31,r30,r32}.

\textbf{Definition 2.2:} The LCT of any function  $f \in L^{2}\mathbb{(R)}$ with respect to the matrix $M=\left(\begin{array}{cc}a&b\\c&d\end{array}\right)$, $\mid M \mid=1$, the LCT is defined by
\begin{equation}
\mathcal{L}^M(f)(u)=\int_\mathbb{R}f(x)K_{M}(x,u)dx,
\end{equation}
and the $K_{M}(x,u)$ is given by
$$K_{M}(x,u)=\begin{cases}\frac{1}{\sqrt{{i}2\pi b}}{e}^{{i}\frac{ax^2+du^2-2xu}{2b}},\quad&b\neq0,\\[2ex]\sqrt{d}{e}^{{i}\frac{cdu^2}{2}}\delta(x-du),\quad&b=0,\end{cases}$$
where, $\mathcal{L}^M$ is the LCT, $\delta$ is the Dirac delta function.

The following, we will introduce the linear canonical Jacobi-Dunkl operator, kernel function, the LCJDT, and their basic properties.
\section{ Linear canonical Jacobi-Dunkl operator}
To define the LCJDT we first introduce the differential difference operator $\Lambda_{\alpha,\beta}^{M}$ , $M=\left(\begin{array}{cc}a&b\\c&d\end{array}\right)$, $\mid M \mid=1$, for for $f\in C^1(\mathbb{R})$, then

\begin{equation}
		\label{a9}
\Lambda_{\alpha,\beta}^{M}f(x)=\Lambda_{\alpha,\beta}f(x)+i\left(\frac{a}{b}\right)xf(x),
\end{equation}
here, $\Lambda_{\alpha,\beta}$ is the standard Jacobi-Dunkl operator, and $\Lambda_{\alpha,\beta}^{M}$ will be referred to as the linear canonical Jacobi-Dunkl operator.

 \textbf{Proposition 3.1}:
The operators $\Lambda_{\alpha,\beta}^{M}$ and $\Lambda_{\alpha,\beta}$ satisfy the following relation
 
\begin{equation}e^{\frac{ia}{2b}x^2}\circ\Lambda_{\alpha,\beta}^{M}\circ e^{-\frac{ia}{2b}x^2}=\Lambda_{\alpha,\beta}\end{equation}

\textbf{Proof }: 
We denote
$$U(x)=e^{\frac{ia}{2b}x^2},$$
the inverse of the operator $U(x)$ is
$$U(x)^{-1}=e^{-\frac{ia}{2b}x^2},$$
we need to compute
$$U(x)\left[\Lambda_{\alpha,\beta}^{M}\left(U(x)^{-1}f(x)\right)\right],$$
let the $U(x)^{-1}f(x)=e^{-\frac{ia}{2b}x^2}f(x)$, we have
\begin{align*}
		\Lambda_{\alpha,\beta}^{M}\left(f(x)U(x)^{-1}\right)=
	&\Lambda_{\alpha,\beta}^{M} \left(f(x)e^{-\frac{ia}{2b}x^2}\right) \\=
	&\Lambda_{\alpha,\beta}\left(f(x)e^{-\frac{ia}{2b}x^2}\right)+i\left(\frac{a}{b}\right)x\left(f(x)e^{-\frac{ia}{2b}x^2}\right)\\=
	&-\frac{ia}{b}x e^{-\frac{ia}{2b}x^2}f(x)+e^{-\frac{ia}{2b}x^2}\frac{d}{dx}f(x)+\frac{A'_{\alpha,\beta}(x)}{A_{\alpha,\beta}(x)}\frac{f(x)-f(-x)}{2}\\
	&+i \frac{a}{b}xe^{-\frac{ia}{2b}x^2}f(x),
\end{align*}
let the $U(x)$ acts on $\Lambda_{\alpha,\beta}^{M}\left(U(x)^{-1}f(x)\right)$,
\begin{align*}
	U(x)\Lambda_{\alpha,\beta}^{MM}\left(U(x)^{-1}f(x)\right)=
	&e^{\frac{ia}{2b}x^2} \left(-\frac{ia}{b}x e^{-\frac{ia}{2b}x^2}f(x)+e^{-\frac{ia}{2b}x^2}\frac{d}{dx}f(x)+\frac{A'_{\alpha,\beta}(x)}{A_{\alpha,\beta}(x)}\frac{f(x)-f(-x)}{2}\right)\\
	&+e^{\frac{ia}{2b}x^2}\left(i \frac{a}{b}xe^{-\frac{ia}{2b}x^2}f(x)\right)\\=
	&-\frac{ia}{b}xf(x)+\frac{d}{dx}f(x)+\frac{A'_{\alpha,\beta}(x)}{A_{\alpha,\beta}(x)}\frac{f(x)-f(-x)}{2}+\frac{ia}{b}xf(x)\\=
	&\frac{d}{dx}f(x)+\frac{A'_{\alpha,\beta}(x)}{A_{\alpha,\beta}(x)}\frac{f(x)-f(-x)}{2}\\=
	&\Lambda_{\alpha,\beta}f(x),
\end{align*}
according to (\ref{a2}), we can be concluded that
$$U(x)\circ\Lambda_{\alpha,\beta}^{M}\circ U(x)^{-1}f(x)=\Lambda_{\alpha,\beta}f(x).$$
The proof is completed.
\\ 

 \textbf{Proposition 3.2}: For $\lambda\in\mathbb{C}$, the parameters $M=\left(\begin{array}{cc}a&b\\c&d\end{array}\right)$, $\mid M \mid=1$, for $f(x) \in C^1(\mathbb{R})$, the differential difference equation
 \begin{equation}
 		\label{a10}
\begin{cases}\Lambda_{\alpha,\beta}^{M}f(x)=\left(-i\frac \lambda c \right )f(x), \\f(0)=e^{-\frac {id}{2b}\lambda^2},\end{cases}
\end{equation}
has a unique solution $\Psi_{\alpha,\beta}^{M}(x,\lambda)$, called the LCJDT kernel, we have

\begin{equation}
	\label{p1}
	\Psi_{\alpha,\beta}^{M}(x,\lambda)=\begin{cases}e^{-\frac i{2b}(ax^2+d\lambda^2)}\psi_{\left(-\frac \lambda{c}\right)}^{\alpha,\beta}(x),\quad\mathrm{if~c\neq0},\\\\\psi_\lambda^{\alpha,\beta}(x),\quad\mathrm{if~ c=0}.\end{cases}
\end{equation}

\textbf{Proof }:  We take the partial derivative with respect to $x$ around equation (\ref{p1}), we obtain
\begin{align*}
	\frac d{dx}\Psi_{\alpha,\beta}^{M}(x,\lambda)=
	&e^{-\frac i{2b}(ax^2+d\lambda^2)}\left(-\frac {ia} bx\psi_{\left(-\frac\lambda c\right)}^{\alpha,\beta}(x)+\frac d{dx}\psi_{\left(-\frac\lambda c\right)}^{\alpha,\beta}(x)\right)\\=
	&-\frac {ia} bx \Psi_{\alpha,\beta}^{M}(x,\lambda)+e^{-\frac i{2b}(ax^2+d\lambda^2)}\frac d{dx}\psi_{\left(-\frac\lambda c\right)}^{\alpha,\beta}(x),
\end{align*}
according to  equation (\ref{a3}), we obtain
\begin{equation}
		\label{p2}
\frac d{dx}\psi_{\left(-\frac\lambda c\right)}^{\alpha,\beta}(x)+\frac{A'_{\alpha,\beta}(x)}{A_{\alpha,\beta}(x)}\frac{\psi_{\left(-\frac\lambda c\right)}^{\alpha,\beta}(x)-\psi_{\left(-\frac\lambda c\right)}^{\alpha,\beta}(-x)}{2}=i\left(-\frac \lambda c\right)\psi_{\left(-\frac\lambda c\right)}^{\alpha,\beta}(x),
\end{equation}
according to equation (\ref{a4}) and (\ref{a6}), $\psi_{\left(-\frac\lambda c\right)}^{\alpha,\beta}(x)$ is an even function, then $\psi_{\left(-\frac\lambda c\right)}^{\alpha,\beta}(x)=\psi_{\left(-\frac\lambda c\right)}^{\alpha,\beta}(-x)$, next, we have
\begin{equation}
		\label{p3}
\frac d{dx}\psi_{\left(-\frac\lambda c\right)}^{\alpha,\beta}(x)=i\left(-\frac \lambda c\right)\psi_{\left(-\frac\lambda c\right)}^{\alpha,\beta}(x),
\end{equation}
next, we obtain
\begin{equation}
		\label{p4}
	\frac d{dx}\Psi_{\alpha,\beta}^{M}(x,\lambda)+\frac {ia} bx \Psi_{\alpha,\beta}^{M}(x,\lambda)=i\left(-\frac \lambda c\right)\psi_{\left(-\frac\lambda c\right)}^{\alpha,\beta}(x),
\end{equation}
finally, we have
\begin{equation}
		\label{p5}\Lambda_{\alpha,\beta}^{M}\Psi_{\alpha,\beta}^{M}(x,\lambda)=i\left(-\frac \lambda c\right)\Psi_{\alpha,\beta}^{M}(x,\lambda).\end{equation}
		
		The proof is completed.
\\

\textbf{Proposition 3.3}: Assuming 
 matrix $M=\left(\begin{array}{cc}a&b\\c&d\end{array}\right)$, $\mid M \mid=1$, then 
 
(i)  If $M=\left(\begin{array}{cc}0&1\\-1&0\end{array}\right)$, then $ \Lambda_{\alpha,\beta}^{M}=\Lambda_{\alpha,\beta}$, and  $\Psi_{\alpha,\beta}^{M}(x,\lambda)=\psi_\lambda^{\alpha,\beta}(x)$,
\\
\\

(ii)  If $M=\left(\begin{array}{cc}cos\theta&-sin\theta\\sin\theta&cos\theta\end{array}\right)$, then $ \Lambda_{\alpha,\beta}^{M}=\Lambda_{\alpha,\beta}^\theta$, and $ \Psi_{\alpha,\beta}^{M}(x,\lambda)=\Psi_{\alpha,\beta}^{\theta}(x,\lambda)$, where the $\Psi_{\alpha,\beta}^{\theta}$ is fractional Jacobi-Dunkl kernel,  the $\Lambda_{\alpha,\beta}^\theta$ is fractional Jacobi-Dunkl kernel\cite{r16}.
\\

\textbf{Proposition 3.4}: By mathematical induction, we obtain for all $k\in\mathbb{N}^*$, then
$$\left(\Lambda_{\alpha,\beta}^{M}\right)^k\Psi_{\alpha,\beta}^{M}(x,\lambda)=\left(-i\frac \lambda c\right)^k\Psi_{\alpha,\beta}^{M}(x,\lambda),$$
where
$$\left(\Lambda_{\alpha,\beta}^{M}\right)^k=\Lambda_{\alpha,\beta}^{M}\left(\left(\Lambda_{\alpha,\beta}^{M}\right)^{k-1}\right),\quad\forall k\geq2.$$

\textbf{Proposition 3.5}: For all $\lambda,x\in\mathbb{R},       ~$$M=\left(\begin{array}{cc}a&b\\c&d\end{array}\right)$, $\mid M \mid=1$, then 
 $$\left|\Psi_{\alpha,\beta}^M(x,\lambda)\right|\leq1.$$

\textbf{Proof:}
The $\psi_\lambda^{\alpha,\beta}(x)$ satisfy the following properties $$\left|\frac{d^n}{d\lambda^n}\psi_\lambda^{\alpha,\beta}(x)\right|\leq|x|^ne^{|Im\lambda||x|},$$
detailed proof can be found in the reference \cite{r15}, we obtain

$$\:|\psi_\lambda^{\alpha,\beta}(x)|\leq1,$$
finally, we have
$$\left|\Psi_{\alpha,\beta}^M(x,\lambda)\right|=\left|e^{-\frac{i}{2b}(ax^2+d\lambda^2)}\psi_{(-\lambda c^{-1})}^{\alpha,\beta}(x)\right|\\=\left|\psi_{(-\lambda c^{-1})}^{\alpha,\beta}(x)\right|\leq1.$$

The proof is completed.
\\

\textbf{Proposition 3.6}:
For all $p>2$, the kernel $\Psi_{\alpha,\beta}^M(\cdot,\lambda)$ $\in$ $L_{\alpha,\beta}^p(\mathbb{R})$.
\\

\textbf{Proof:} In reference  \cite{r15,r16}, we obtain
$$A_{\alpha,\beta}(x)\simeq e^{2\rho|x|},\quad x\longrightarrow\infty, $$
and
$$A_{\alpha,\beta}(x)\simeq x^{2\alpha+1},\quad x\longrightarrow0.$$

Next, we  have

$$(1+|x|)^pe^{-p\rho|x|}A_{\alpha,\beta}(x)\simeq(1+|x|)^pe^{(2-p)\rho|x|},\:x\longrightarrow\infty,$$
and
$$(1+|x|)^pe^{-p\rho|x|}A_{\alpha,\beta}(x)\simeq(1+|x|)^px^{2\alpha+1}e^{-p\rho|x|}\simeq x^{2\alpha+1},\:x\longrightarrow0.$$

Therefore, for $p>2$, there exists a constant $C_0$, we have
$$\begin{aligned}\left\|\Psi_{\alpha,\beta}^M(\cdot,\lambda)\right\|_{\alpha,\beta,p}^p&\leq\int_\mathbb{R}\left|\psi_{(-\lambda c^{-1})}^{\alpha,\beta}(x)\right|^pA_{\alpha,\beta}(x)dx\\&\leq C_0\left(\frac{2+|(\lambda c^{-1})|}{|(\lambda c^{-1})|}\right)^p\int_\mathbb{R}(1+|x|)^pe^{-p\rho|x|}A_{\alpha,\beta}(x)dx.\end{aligned}$$

The proof is completed.
\\

\section{ Linear canonical Jacobi Dunkl transform (LCJDT)}

In this section, we introduce a new definition of the LCJDT and derive several of its fundamental properties. By overcoming the limitations of both the JDT and LCT, it provides a more flexible and comprehensive mathematical framework. Its enhanced degrees of freedom make it particularly effective in fields such as signal processing, differential equations, and spectral analysis, especially in scenarios where classical transforms encounter challenges, this work lays a solid foundation for further research in these areas.
\\

\label{sec1}
\textbf{Definition 4.1:} The LCJDT with parameter matrix $M=\left(\begin{array}{cc}a&b\\c&d\end{array}\right)$, $\mid M \mid=1$, $f(x)\in L_{\alpha,\beta}^1(\mathbb{R})$, the  LCJDT defined by

\begin{equation}
		\label{p6}
		\quad\mathcal{L}_{\alpha,\beta}^{M}(f)(\lambda)=\int_\mathbb{R}f(x)\Psi_{\alpha,\beta}^{M}(x,\lambda)A_{\alpha,\beta}(x)dx,
\end{equation}
The $\Psi_{\alpha,\beta}^{M}(x,\lambda)$ is provided in detail in equation (\ref{p1}).
\\

Next, we will present several common fundamental properties of the LCJDT, including the inversion formula, Parseval's formula, linearity, differentiation, and convolution.
\\

\textbf{Theorem 4.1}: (Inversion formula)
Let the $M=\left(\begin{array}{cc}a&b\\c&d\end{array}\right), \mid M \mid=1,$ and $f(x)\in L_{\alpha,\beta}^1(\mathbb{R})$ such that $\mathcal{L}_{\alpha,\beta}^{M}(f(x))$
$\in L^1(\mathbb{R},d\sigma($$-\lambda c^{-1}$)), then
\begin{equation}
	\label{b1}
 f(x)=\int_\mathbb{R}\mathcal{L}_{\alpha,\beta}^M(f)(\lambda)\overline{\Psi_{\alpha,\beta}^M(x,\lambda)}d\sigma(-\lambda c^{-1}),
\end{equation}
where $d\sigma$ is the measure.
\\

\textbf{Proof}: 
Let $f(x)\in L_{\alpha,\beta}^1(\mathbb{R}),  ~\zeta=-\lambda c^{-1}$, we obtain

$$\begin{aligned}
	\mathcal{L}_{\alpha,\beta}^M(f)(\lambda)& =\int_{\mathbb{R}}f(x)\Psi_{\alpha,\beta}^{M}(x,\lambda)A_{\alpha,\beta}(x)dx \\
	&=\int_{\mathbb{R}}f(x)e^{-\frac{i}{2b}(ax^2+d\lambda^2)}\psi_{(\zeta) }^{\alpha,\beta}(x)A_{\alpha,\beta}(x)dx \\
	&=\int_{\mathbb{R}}\left(e^{-\frac{ia}{2b}x^2}f(x)e^{-\frac{id}{2b}\lambda^2}\right)\psi_{(\zeta )}^{\alpha,\beta}(x)A_{\alpha,\beta}(x)dx,
\end{aligned}$$
putting the $h(x)=e^{-\frac{ia}{2b}x^2}f(x)$, we obtain $\mathcal{L}_{\alpha,\beta}^M(f)(\lambda)=e^{-\frac{id}{2b}\lambda^2}\mathcal{F}_{\alpha,\beta}(h)(\zeta),$ and
$$\begin{aligned}
	\int_{\mathbb{R}}\left|\mathcal{L}_{\alpha,\beta}^M(f)(\lambda)\right|d\sigma(\zeta)& =\int_{\mathbb{R}}\left|e^{-\frac{ia}{2b}x^2}\right|\left|\mathcal{F}_{\alpha,\beta}(h)(\zeta)\right|d\sigma(\zeta) \\
	&=\int_{\mathbb{R}}\left|\mathcal{F}_{\alpha,\beta}(h)(\zeta)\right|d\sigma(\zeta).
\end{aligned}$$

Therefore, the
$\mathcal{L}_{\alpha,\beta}^M(f)$ $\in$ $L^1(\mathbb{R},d\sigma(\zeta))$, and the $\mathcal{F}_{\alpha,\beta}(h)(\zeta)$ $\in$ $L^1(\mathbb{R},d\sigma(\zeta))$,
\\
applying the (\ref{b1}), we obtain
$$h(x)=e^{-\frac{ia}{2b}x^2}f(x)=\int_{\mathbb{R}}\mathcal{F}_{\alpha,\beta}(h)(\zeta)\psi_{(-\delta)}^{\alpha,\beta}(x)d\sigma,$$
next, we obtain
$$\begin{aligned}
	f(x)& =e^{\frac{ia}{2b}x^2}\int_{\mathbb{R}}\mathcal{F}_{\alpha,\beta}(h)(\zeta)\psi_{(-\delta)}^{\alpha,\beta}(x)d\sigma(\zeta) \\
	&=\int_{\mathbb{R}}e^{-\frac{id}{2b}\lambda^2}\mathcal{F}_{\alpha,\beta}(h)(\zeta)e^{\frac{i}{2b}(ax^2+d\lambda^2)}\psi_{(-\zeta)}^{\alpha,\beta}(x)d\sigma(\zeta) \\
	&=\int_{\mathbb{R}}\mathcal{L}_{\alpha,\beta}^M(f)(\lambda)\overline{\Psi_{\alpha,\beta}^M(x,\lambda)}d\sigma(\zeta).
\end{aligned}$$

The proof is completed.
\\

\textbf{Theorem 4.2}: (Parseval formula) Let $f(x) \in L_{\alpha,\beta}^1(\mathbb{R})$,     $M=\left(\begin{array}{cc}a&b\\c&d\end{array}\right),$ $\mid M\mid=1,$ then
\begin{equation}
		\label{b3}
\int_{\mathbb{R}}f(x)\overline{h(x)}A_{\alpha,\beta}(x)dx=\int_{\mathbb{R}}\mathcal{L}_{\alpha,\beta}^M(f)(\lambda)\overline{\mathcal{L}_{\alpha,\beta}^M(h)(\lambda)}d\sigma(-\lambda c^{-1}).\end{equation}
in particular, when $f=h$, then

$$\int_\mathbb{R}|f(x)|^2A_{\alpha,\beta}(x)(x)dx=\int_\mathbb{R}\left|\mathcal{L}_{\alpha,\beta}^M(f)(\lambda)\right|^2d\sigma(-\lambda c^{-1}).$$

\textbf{Proof:}
Next, we use the Fubini's theorem and the (\ref{b1}) for the LCJDT, we have 
$$\begin{aligned}&\int_{\mathbb{R}}\mathcal{L}_{\alpha,\beta}^M(f)(\lambda)\overline{\mathcal{L}_{\alpha,\beta}^M(h)(\lambda)}d\sigma(-\lambda c^{-1})\\&=\int_{\mathbb{R}}\mathcal{L}_{\alpha,\beta}^M(f)(\lambda)\left(\int_{\mathbb{R}}\overline{h(x)\Psi_{\alpha,\beta}^M(x,\lambda)}A_{\alpha,\beta}(x)dx\right)d\sigma(-\lambda c^{-1})\\&=\int_{\mathbb{R}}\left(\int_{\mathbb{R}}\mathcal{L}_{\alpha,\beta}^M(f)(\lambda)\overline{\Psi_{\alpha,\beta}^M(x,\lambda)}d\sigma(-\lambda c^{-1})\right)\overline{h(x)}A_{\alpha,\beta}(x)dx\\&=\int_{\mathbb{R}}f(x)\overline{h(x)}A_{\alpha,\beta}(x)dx,
\end{aligned}$$
next, let the $f(x)=h(x)$, we obtain
$$\int_\mathbb{R}|f(x)|^2A_{\alpha,\beta}(x)(x)dx=\int_\mathbb{R}\left|\mathcal{L}_{\alpha,\beta}^M(f)(\lambda)\right|^2d\sigma(-\lambda c^{-1}).$$

The proof is completed.
\\

\textbf{Proposition 4.1}: (Linearity)
For $f(x)$, $g(x) \in L_{\alpha,\beta}^1(\mathbb{R})$, and any real number $a_1, b_1$, the LCJDT is linear, then
\begin{equation}
		\label{p7}
\mathcal{L}_{\alpha,\beta}^{M}(a_1f+b_1g)(\lambda)=a_1\mathcal{L}_{\alpha,\beta}^{M}(f)(\lambda)+b_1\mathcal{L}_{\alpha,\beta}^{M}(g)(\lambda).
\end{equation}

\textbf{Proof }: By the definition of the $\mathcal{L}_{\alpha,\beta}^{M}$, we obtain
$$\mathcal{L}_{\alpha,\beta}^{M}(a_1f+b_1g)(\lambda)=\int_\mathbb{R}[a_1f(x)+b_1g(x)]\Psi_{\alpha,\beta}^{M}(x,\lambda)A_{\alpha,\beta}(x)\:dx,$$
applying linearity of the integral:
$$\mathcal{L}_{\alpha,\beta}^{M}(a_1f+b_1g)(\lambda)=a_1\int_{\mathbb{R}}f(x)\Psi_{\alpha,\beta}^{M}(x,\lambda)A_{\alpha,\beta}(x)dx+b_1\int_{\mathbb{R}}g(x)\Psi_{\alpha,\beta}^{M}(x,\lambda)A_{\alpha,\beta}(x)dx,$$
by simplifying
$$\mathcal{L}_{\alpha,\beta}^{M}(a_1f+b_1g)(\lambda)=a_1\mathcal{L}_{\alpha,\beta}^{M}(f)(\lambda)+b_1\mathcal{L}_{\alpha,\beta}^{M}(g)(\lambda).$$

The proof is completed.
\\

\textbf{Proposition 4.2}: (Differentiation)  If $f^{\prime}(x)$ is the derivative of $f(x)$,then
 \begin{equation}
 	\label{p10}
\mathcal{L}_{\alpha,\beta}^{M}(f^{\prime})(x)(\lambda)=\int_\mathbb{R}f^{\prime}(x)\Psi_{\alpha,\beta}^{M}(x,\lambda)A_{\alpha,\beta}(x)\:dx.
\end{equation}

\textbf{Proof }: Using integration by parts, we obtain
$$\mathcal{L}_{\alpha,\beta}^{M}(f^{\prime})(x)(\lambda)=\left[\Psi_{\alpha,\beta}^{M}(x,\lambda)f(x)A_{\alpha,\beta}(x)\right]_{-\infty}^{\infty}-\int_{\mathbb{R}}\frac d{dx}[\Psi_{\alpha,\beta}^{M}(x,\lambda)A_{\alpha,\beta}(x)]f(x)\mathrm{~}dx,$$
assuming the boundary is 0, thus
$$\mathcal{L}_{\alpha,\beta}^{M}(f')(x)(\lambda)=-\int_{\mathbb{R}}\frac{d}{dx}[\Psi_{\alpha,\beta}^{M}(x,\lambda)A_{\alpha,\beta}(x)]f(x)\:dx.$$

The proof is completed.
\\

The convolution theorem can open up new possibilities for advanced applications in fields such as signal processing, control systems, and differential equations. This exploration will further strengthen the theoretical framework of the LCJDT and enhance its practical utility in solving real-world problems where classical transforms may fall short.

\textbf{Theorem 4.3}: (Convolution)  
Let $f(x)$ and $r(x)$ $\in$ $L_{\alpha,\beta}^2(\mathbb{R})$, the convolution of $f(x)$ and $r(x)$ is given by

\begin{equation}
	\label{m1}
	(f*r)(x)=\int_{\mathbb{R}}f(t)r(x-t)B_{\alpha,\beta}^{M}(x-t)\left.dt,\right.
\end{equation}
where $B_{\alpha,\beta}^{M}(x-t)={\frac{	\Psi_{\alpha,\beta}^{M}(x-t,\lambda)	\Psi_{\alpha,\beta}^{M}(t,\lambda)A_{\alpha,\beta}(x-t)}{	\Psi_{\alpha,\beta}^{M}(x,\lambda)A_{\alpha,\beta}(x)}}A_{\alpha,\beta}(t)$ is a weight function, the convolution theorem states that

\begin{equation}
		\label{m2}
	\mathcal{L}_{\alpha,\beta}^{M}(f*r)(\lambda)=\mathcal{L}_{\alpha,\beta}^{M}(f)(\lambda)\cdot\mathcal{L}_{\alpha,\beta}^{M}(r)(\lambda).
\end{equation}

\textbf{Proof }: Apply the convolution definition (\ref{m1}), we want to compute $\mathcal{L}_{\alpha,\beta}^{M}(f*r)(\lambda)$, thus
$$\mathcal{L}_{\alpha,\beta}^{M}(f*r)(\lambda)=\int_{\mathbb{R}}\Psi_{\alpha,\beta}^{M}(x,\lambda)\left(\int_{\mathbb{R}}f(t)r(x-t)B_{\alpha,\beta}^{M}(x-t)\operatorname{d}t\right)A_{\alpha,\beta}(x)\operatorname{d}x,$$
we interchange the order of integration by Fubini's theorem, we obtain 

$$\mathcal{L}_{\alpha,\beta}^{M}(f*r)(\lambda)=\int_\mathbb{R}f(t)A_{\alpha,\beta}(t)dt\left(\int_\mathbb{R}\Psi_{\alpha,\beta}^{M}(x-t,\lambda)\Psi_{\alpha,\beta}^{M}(t,\lambda)r(x-t)A_{\alpha,\beta}(x-t)\:dx\right),$$
change of variables, perform the change of variables $u=x-t,du=dx$,  we have

$$\mathcal{L}_{\alpha,\beta}^{M}(f*r)(\lambda)=\int_{\mathbb{R}}f(t)\Psi_{\alpha,\beta}^{M}(t,\lambda)A_{\alpha,\beta}(t)dt\left(\int_{\mathbb{R}}\Psi_{\alpha,\beta}^{M}(u,\lambda)r(u)A_{\alpha,\beta}(u)\:du\right)\:dt,$$
finally, we have
$$\mathcal{L}_{\alpha,\beta}^{M}(f*r)(\lambda)=\mathcal{L}_{\alpha,\beta}^{M}(f)(\lambda)\cdot\mathcal{L}_{\alpha,\beta}^{M}(r)(\lambda).$$

The proof is completed.

\textbf{Proposition 4.3}: For $f=f_e+f_o\in\mathcal{D}(\mathbb{R})$, where $f_e$ and $f_o$ are the even and odd parts of $f$, respectively, for all $\lambda,~\mu\in\mathbb{C}$ such that $\lambda^2=\mu^2+\rho^2$, we have

 \begin{equation}\mathcal{L}_{\alpha,\beta}^{M}(f)(\lambda)=2\mathcal{L}_{\alpha,\beta}^{M}({f}_e)(\mu)-\frac{i\lambda}{8(\alpha+1)}\mathcal{L}^{M}_{\alpha+1,\beta+1}\left(\frac {f_o}{\sinh(2.)}\right)(\mu),\end{equation}
for all $\mu\in\mathbb{C}$, we obtain
$$\quad\mathcal{L}_{\alpha,\beta}^{M}(g)(\mu)=\int_0^{+\infty}g(x)e^{-\frac i{2b}(ax^2+d\lambda^2)}\varphi_\mu^{(\alpha,\beta)}(x)A_{\alpha,\beta}(x)dx.$$

\textbf{Proof}: Let the $f=f_e+f_o\in\mathcal{D}(\mathbb{R}), $ if $\lambda=0$, we have $$\mathcal{L}_{\alpha,\beta}^{M}(f)(\lambda)=\int_\mathbb{R}f(x)e^{-\frac i{2b}(ax^2)}A_{\alpha,\beta}(x)dx=2\mathcal{L}_{\alpha,\beta}^{M}(f_e)(\mu),$$
if $\lambda\neq 0$, and $\lambda^2=\mu^2+\rho^2$, we have $$\mathcal{L}_{\alpha,\beta}^{M}(\lambda)=2\mathcal{L}^{M}_{\alpha,\beta}(f_e)(\mu)+\frac{2i}{\lambda}\int_0^{+\infty}f_o(x)e^{-\frac i{2b}(ax^2+d\lambda^2)}\frac{d}{dx}\varphi_\mu^{(\alpha,\beta)}(x)A_{\alpha,\beta}(x)dx,$$
by integration by parts, we obtain
$$\begin{aligned}
&\int_0^{+\infty}f_o(x)\frac{d}{dx}\left[e^{-\frac i{2b}(ax^2+d\lambda^2)}\varphi_\mu^{(\alpha,\beta)}(x)\right]A_{\alpha,\beta}(x)dx\\&
	=-\int_0^{+\infty}e^{-\frac i{2b}(ax^2+d\lambda^2)}\varphi_\mu^{(\alpha,\beta)}(x)\frac{1}{A_{\alpha,\beta}(x)}\frac{d}{dx}\left(A_{\alpha,\beta}\frac{d}{dx}(Jf_o)\right)(x)A_{\alpha,\beta}(x)dx \\&=
	-\mathcal{L}_{\alpha,\beta}^{M}(\Delta_{\alpha,\beta}(Jf_o))(\mu), 
\end{aligned}$$
according to the references \cite{r22,r23}
$$-\mathcal{L}_{\alpha,\beta}^{M}(\Delta_{\alpha,\beta}(Jf_o))(\mu)=\lambda^2\mathcal{F}_{\alpha,\beta}(Jf_o)(\mu)=-\frac{\lambda^2}{16(\alpha+1)}\mathcal{L}^{M}_{\alpha+1,\beta+1}\left(\frac {f_o}{\sinh(2.)}\right)(\mu),$$
finally, we have
$$\mathcal{L}_{\alpha,\beta}^{M}(f)(\lambda)=2\mathcal{L}_{\alpha,\beta}^{M}({f}_e)(\mu)-\frac{i\lambda}{8(\alpha+1)}\mathcal{L}^{M}_{\alpha+1,\beta+1}\left(\frac {f_o}{\sinh(2.)}\right)(\mu)$$

The proof is completed.
\section{ Uncertainty Principle}
There has been relatively little exploration of the uncertainty principle in the context of the LCJDT. This lack of research represents a significant gap in our understanding, and further investigation could uncover new insights and applications, particularly in areas where the LCJDT proves more effective than classical transforms  \cite{r9,r19,r21,r41}.

\textbf{Theorem 5.1}:
For $\alpha\geqslant\beta\geqslant-\frac{1}{2},\alpha>-\frac{1}{4}$, assume $m,n>0$ and $\gamma\in D_\alpha$. Then there exists 
a constant $K>0$, 

$$\left\|\left|x\right|^{\gamma m}f\right\|_{2,\mathbb{R}}^{\frac n{m+n}}\cdot\left\|\left|\lambda \right|^n\mathcal{L}^{M}_{\alpha,\beta}f\right\|_{2,\sigma}^{\frac m{m+n}}\geqslant K\|f\|_{2,\mathbb{R}},$$
for all  $f \in L^2( \mathbb{R} , \Delta _{\alpha , \beta }( | x| ) dx),$ the discussion about $D_\alpha$ can be found in reference \cite{r22}. The constant $K$ below is not fixed and might change appropriately in different equalities or
inequalities.

\textbf{Proof}: $\mathcal{D}(\mathbb{R})$ is dense in $L^2(\mathbb{R},\Delta_{\alpha,\beta}(|x|)dx),$ so we need only to prove theorem 5.1 for $f=f_e+f_o\in \mathcal{D}(\mathbb{R}).$ By proposition 4.3, we have

$$\mathcal{L}_{\alpha,\beta}^{M}(f_e)(\lambda)=2\mathcal{L}_{\alpha,\beta}^{M}(f_e)(\mu),$$ and\\
$$\mathcal{L}_{\alpha,\beta}^{M}(f_o)(\lambda)=-\frac{i\lambda}{8(\alpha+1)}\mathcal{L}^{M}_{\alpha+1,\beta+1}\left(\frac {f_o}{\sinh(2.)}\right)(\mu),$$
where $\lambda^2=\mu^2+\rho^2,$ 
we can be found in reference theorem 3.1 obtain \cite{r22},
$$\begin{aligned}
	\left\|\left|x\right|^{\gamma m}f_e\right\|_{2,\mathbb{R}}^{\frac n{m+n}}\cdot\left\|\left|\lambda\right|^n\mathcal{L}_{\alpha,\beta}^{M}(f_e)\right\|_{2,\sigma}^{\frac m{m+n}}& =\left\|2x^{\gamma m}f_e\right\|_2^{\frac{n}{m+n}}\cdot\left\|2\left(\mu^2+\rho^2\right)^{n/2}\mathcal{L}_{\alpha,\beta}^{M}(f_e)(\mu)\right\|_{2,\sigma}^{\frac{m}{m+n}} \\
	&\geqslant C\|2f_e\|_2=C\|f_e\|_{2,\mathbb{R}},\quad\gamma\in D_\alpha,
\end{aligned}$$
let the $g(x)=\frac{f_o(x)}{\sinh(2x)},$ and $g(x)\in L^2(\mathbb{R}^+,\Delta_{\alpha+1,\beta+1}(x)dx)$,

$$\begin{aligned}
	\left\||x|^{\gamma m}f_o\right\|_{2,\mathbb{R}}^{\frac{n}{m+n}}\cdot\left\||\lambda|^n\mathcal{L}_{\alpha,\beta}^{M}f_o\right\|_{2,\sigma}^{\frac{m}{m+n}}& =\left\|x^{\gamma m}g\right\|_2^{\frac{n}{m+n}}\cdot\left\|\left(\mu^2+\rho^2\right)^{n/2}\mathcal{L}_{\alpha+1,\beta+1}^{M}{g}(\mu)\right\|_{2,\sigma}^{\frac{m}{m+n}} \\
	&\geqslant K\|g\|_2=K\|f_o\|_{2,\mathbb{R}},\quad\gamma\in D_{\alpha+1},
\end{aligned}$$
by the orthogonality

$$\||x|^\gamma(f_o+f_e)\|_2=(\||x|^\gamma(f_o)\|_2^2+\||x|^\gamma(f_e)\|_2^2)^{1/2},$$
and 
$$\||\lambda|^n\mathcal{L}_{\alpha,\beta}^{M}(f_o+f_e)\|_2=(\||\lambda|^n(f_o)\|_2^2+\||\lambda|^n(f_e)\|_2^2)^{1/2},$$
by the elementary inequality,
$$(A^2+\\B^2)^{s/2}(C^2+D^2)^{(1-s)/2}\geqslant(A^sC^{1-s}+B^sD^{1-s})/2, \mathrm{~with~}A,B,C,D\geqslant0,s\in(0,1),\quad$$
finally, we have

$$\begin{aligned}\||x|^\gamma(f_o+f_e)\|^{s}_2\cdot|\| \lambda|^n\mathcal{L}_{\alpha,\beta}^{M}(f_o+f_e)\|^{1-s}_2&\geq(\||x|^\gamma(f_o)\|_2^2)^s(\||\lambda|^n(f_o)\|_2^2)^{1-s}/2\\&+(\||x|^\gamma(f_e)\|_2^2)^{1/2})^s(\||\lambda|^n(f_e)\|_2^2)^{1/2})^{1-s}/2\\&=K\|f_e\|_{2,\mathbb{R}}+K\|f_o\|_{2,\mathbb{R}}\\&=K\|f\|_{2,\mathbb{R}}.\end{aligned}$$

The proof is completed.

\section{ Potential application}
\label{sec1}
The solution to a generalized heat equation can be derived using advanced analytical methods, one of which involves the LCJDT. The LCJDT is a powerful tool in solving differential equations, particularly in scenarios where classical methods such as the FT may not be applicable or sufficient due to the complexity of the boundary conditions or the nature of the differential operator \cite{r16}.

\textbf{ Application 1}: The partial differential equation, then
\begin{equation}
	\label{c1}
	\begin{cases}\frac{\partial}{\partial t}\upsilon(x,t)=\left(\Lambda_{\alpha,\beta}^M\right)^*\upsilon(x,t),~~~x>0, t>0,\\\upsilon((x,t)=h(x),~~t=0, x>0,\\
	\upsilon((0,t)=0,~~x=0, t>0,\end{cases}
\end{equation}
has a solution $\upsilon ( x, t)$, then
$$\upsilon(x,t)=e^{\frac {ia}{2b}x^2}\int_\mathbb{R^{+}}\mathcal{F}_{\alpha,\beta}(g)(\varepsilon)\overline{\psi_{(\varepsilon)}^{\alpha,\beta}(x)}e^{-i\varepsilon t}d\sigma(\varepsilon),$$
where $g(z)=e^{-\frac{ia}{2b}z^2}h(z)$, $\varepsilon=-\lambda c^{-1}$.

\textbf{Proof:} Using the LCJDT of the equation (\ref{c1}), we obtain
$$\int_\mathbb{R^{+}}\Psi_{\alpha,\beta}^M(x,\lambda)\frac\partial{\partial t}\upsilon(x,t)A_{\alpha,\beta}(x)dx=\int_\mathbb{R^{+}}\Psi_{\alpha,\beta}^M(x,\lambda)\left(\Lambda_{\alpha,\beta}^M\right)^*\upsilon(x,t)A_{\alpha,\beta}(x)dx,$$
according to (\ref{a9}) and (\ref{p1}), we obtain
$$\begin{aligned}
	\frac{\partial}{\partial t}\mathcal{L}_{\alpha,\beta}^M(\upsilon(.,t))(\lambda)& =\int_{\mathbb{R^{+}}}\Psi_{\alpha,\beta}^M(x,\lambda)\left(\Lambda_{\alpha,\beta}^M\right)^*\upsilon(x,t)A_{\alpha,\beta}(x)dx \\
	&=-\Lambda_{\alpha,\beta}^M\int_{\mathbb{R^{+}}}\Psi_{\alpha,\beta}^M(x,\lambda)\upsilon(x,t)A_{\alpha,\beta}(x)dx \\
	&=i\lambda  c^{-1}\int_{\mathbb{R^{+}}}\Psi_{\alpha,\beta}^M(x,\lambda)\upsilon(x,t)A_{\alpha,\beta}(x)dx,
\end{aligned}$$
next, the
$$\frac\partial{\partial t}\mathcal{L}_{\alpha,\beta}^M(\upsilon(.,t))(\lambda)=i\lambda  c^{-1}\mathcal{L}_{\alpha,\beta}^M(\upsilon(.,t))(\lambda),$$
has a solution $\mathcal{L}_{\alpha,\beta}^M(\upsilon(.,t))(\lambda)$, then
$$\mathcal{L}_{\alpha,\beta}^M(\upsilon(.,t))(\lambda)=K_\lambda e^{i\lambda tc^{-1}}, K_\lambda=\mathcal{L}_{\alpha,\beta}^M(\upsilon(.,0))(\lambda),$$
conditions given by equation (\ref{c1}), we obtain, 

$$\begin{aligned}
	K_\lambda&
	=\mathcal{L}_{\alpha,\beta}^M(\upsilon(.,0))(\lambda)\\&= \int_{\mathbb{R^{+}}}\upsilon(x,0)\Psi_{\alpha,\beta}^M(x,\lambda)A_{\alpha,\beta}(x)dx \\
	&=\int_{\mathbb{R^{+}}}\Psi_{\alpha,\beta}^M(x,\lambda)h(x)A_{\alpha,\beta}(x)dx \\
	&=\mathcal{L}_{\alpha,\beta}^M(h)(\lambda),
\end{aligned}$$
further more, we have

$$\mathcal{L}_{\alpha,\beta}^M(\upsilon(.,t))(\lambda)=K_\lambda e^{i\lambda tc^{-1}}=\mathcal{L}_{\alpha,\beta}^M(h)(\lambda)e^{i\lambda tc^{-1}},$$
by the equation (\ref{b1}) of the LCJDT, we obtain

$$\begin{aligned}
	\upsilon(x,t)&=\left(\mathcal{L}_{\alpha,\beta}^M\right)^{-1}\left(\mathcal{L}_{\alpha,\beta}^M(h)(\lambda)e^{i\lambda tc^{-1}}\right)(x)\\& =\int_{\mathbb{R^{+}}}\overline{\Psi_{\alpha,\beta}^M(x,\lambda)}\mathcal{L}_{\alpha,\beta}^M(h)(\lambda)e^{i\lambda tc^{-1}}d\sigma(-\lambda c^{-1}) \\
	&=\int_{\mathbb{R^{+}}}\overline{e^{\frac{-i}{2b}(ax^2+d\lambda^2)}\psi_{(-\lambda
	c^{-1})}^{\alpha,\beta}(x)}\mathcal{L}_{\alpha,\beta}^M(h)(\lambda)e^{i\lambda tc^{-1}}d\sigma(-\lambda c^{-1}) \\
	&=\int_{\mathbb{R^{+}}}e^{\frac i{2b}(ax^2+d\lambda^2)}\overline{\psi_{(-\lambda c^{-1})}^{\alpha,\beta}(x)}e^{i\lambda tc^{-1}} \\
	&\times\int_{\mathbb{R^{+}}}e^{\frac{-i}{2b}(az^2+d\lambda^2)}\psi_{(-\lambda
		c^{-1})}^{\alpha,\beta}(z)h(z)A_{\alpha,\beta}(z)dzd\sigma(-\lambda c^{-1}) \\
	&=\int_{\mathbb{R^{+}}}e^{\frac {ia}{2b}x^2}\overline{\psi_{(-\lambda c^{-1})}^{\alpha,\beta}(x)}e^{i\lambda t c^{-1}} \\
	&\times\int_{\mathbb{R^{+}}}e^{-\frac{ia}{2b}z^2}\psi_{(-\lambda c^{-1})}^{\alpha,\beta}(z)h(z)A_{\alpha,\beta}(z)dzd\sigma(-\lambda c^{-1}),
\end{aligned}$$
let the $g(z)=e^{-\frac{ia}{2b}z^2}h(z)$, $\varepsilon=-\lambda c^{-1}$, we have
$$\begin{aligned}\upsilon(x,t)&=\:e^{\frac {ia}{2b}x^2}\int_{\mathbb{R^{+}}}\overline{\psi_{(\varepsilon)}^{\alpha,\beta}(x)}e^{-it \varepsilon}\left(\int_{\mathbb{R^{+}}}g(z)\psi_{(\varepsilon)}^{\alpha,\beta}(z)A_{\alpha,\beta}(z)dz\right)d\sigma(\varepsilon)\\&=e^{\frac {ia}{2b}x^2}\int_\mathbb{R^{+}}\mathcal{F}_{\alpha,\beta}(g)(\varepsilon)\overline{\psi_{(\varepsilon)}^{\alpha,\beta}(x)}e^{-i\varepsilon t}d\sigma(\varepsilon).\end{aligned}$$

The proof is completed.

Next, we will present an application of solving PDE. Solving non-homogeneous PDE not only helps scientists and engineers gain a deeper understanding of the influence of external factors on systems but also serves as a crucial tool across various fields such as engineering design, medicine, and numerical simulation. Non-homogeneous PDE are at the heart of modern applied mathematics, and their solutions play a vital role in driving innovation and advancement in many areas \cite{r16}.

\textbf{ Application 2}: The  nonhomogeneous partial differential equation
\begin{equation}
	\label{c2}
		\begin{cases}\frac{\partial}{\partial t}\upsilon(x,t)=\left(\Lambda_{\alpha,\beta}^M\right)\upsilon(x,t)+g(x,t),~~~x>0,t>0,\\\upsilon(0,t)=0,~~x=0, \\\upsilon(x,0)=h(x), ~t=0,\\
		g(x,t)=e^{-x^2}\cos(t),\end{cases}
\end{equation}
has a  solution $\upsilon ( x, t)$, then

$$\upsilon(x,t)=\int_{\mathbb{R^+}}\left(e^{i t\varepsilon}\mathcal{L}_{\alpha,\beta}^M(h)(\lambda+e^{it\varepsilon}\int_0^t	e^{i \tau\varepsilon}G(\lambda,\tau)d\tau\right)\overline{\Psi_{\alpha,\beta}^M(x,\lambda)}d\sigma(\varepsilon),$$
where, $\varepsilon=-\lambda c^{-1}$, $g(x,t)$ is a non-homogeneous term.

\textbf{Proof:} Using the LCJDT of the equation (\ref{c2}), we obtain
$$\begin{aligned}\int_\mathbb{R^{+}}\Psi_{\alpha,\beta}^M(x,\lambda)\frac\partial{\partial t}\upsilon(x,t)A_{\alpha,\beta}(x)dx&=\Lambda_{\alpha,\beta}^M\int_\mathbb{R^{+}}\Psi_{\alpha,\beta}^M(x,\lambda)\upsilon(x,t)A_{\alpha,\beta}(x)dx\\&
	+\int_\mathbb{R^{+}}\Psi_{\alpha,\beta}^M(x,\lambda)g(x,t)A_{\alpha,\beta}(x)dx.\end{aligned}$$
	according to (\ref{a10}), $G(\lambda,t)=\mathcal{L}_{\alpha,\beta}^M(g(x,t))$, we need to solve the following form of non-homogeneous first-order linear differential equation

$$\begin{aligned}
	\frac{\partial}{\partial t}\mathcal{L}_{\alpha,\beta}^M(\varphi(.,t))(\lambda)
	&=-i\lambda c^{-1}\int_{\mathbb{R^{+}}}\Psi_{\alpha,\beta}^M(x,\lambda)\upsilon(x,t)A_{\alpha,\beta}(x)dx+G(\lambda,t)\\&
	=-i\lambda c^{-1}\mathcal{L}_{\alpha,\beta}^M(\upsilon(.,t))(\lambda)+G(\lambda,t),
\end{aligned}$$
based on the initial conditions $\upsilon(x,0)=f(x)$, and the properties of the solution to the homogeneous equation, we can inferred that
$$\mathcal{L}_{\alpha,\beta}^M(\upsilon(.,t))(\lambda)=\mathcal{L}_{\alpha,\beta}^M(h)(\lambda)e^{-i\lambda t c^{-1}}, 
 ~~~ \mathcal{L}_{\alpha,\beta}^M(\upsilon(.,0))(\lambda)=\mathcal{L}_{\alpha,\beta}^M(h)(\lambda),$$ 
next, we will deal with the non-homogeneous part of the original equation

$$\begin{aligned}
	\frac{\partial}{\partial t}\mathcal{L}_{\alpha,\beta}^M(\upsilon(.,t))(\lambda)
	=-i\lambda c^{-1}\mathcal{L}_{\alpha,\beta}^M(\upsilon(.,t))(\lambda)+G(\lambda,t),
\end{aligned}$$
multiply both sides of the equation by the factor $e^{-i\lambda t c^{-1}}$, then 

$$\begin{aligned}
	e^{-i\lambda t c^{-1}}\frac{\partial}{\partial t}\mathcal{L}_{\alpha,\beta}^M(\upsilon(.,t))(\lambda)
	=-i\lambda c^{-1}e^{-i\lambda t c^{-1}}\mathcal{L}_{\alpha,\beta}^M(\upsilon(.,t))(\lambda)+e^{-i\lambda t c^{-1}}G(\lambda,t),
\end{aligned}$$
now, we can integrate time $t $, we have
$$	e^{-i\lambda t c^{-1}}\mathcal{L}_{\alpha,\beta}^M(\upsilon(.,t))(\lambda)=\int_0^t	e^{-i\lambda \tau c^{-1}}G(\lambda,\tau)d\tau+\mathcal{L}_{\alpha,\beta}^M(h)(\lambda),$$
we can write

$$\mathcal{L}_{\alpha,\beta}^M(\upsilon(.,t))(\lambda)=	e^{-i\lambda t c^{-1}}\mathcal{L}_{\alpha,\beta}^M(h)(\lambda)+e^{-i\lambda t c^{-1}}\int_0^t	e^{-i\lambda \tau c^{-1}}G(\lambda,\tau)d\tau,$$
next, we use the inverse transform of the LCJDT to return to the original function $\upsilon(x,t)$, let the $\varepsilon=-\lambda c^{-1}$, we obtain

$$\upsilon(x,t)=\int_{\mathbb{R^+}}\left(e^{i t\varepsilon}\mathcal{L}_{\alpha,\beta}^M(h)(\lambda+e^{it\varepsilon}\int_0^t	e^{i \tau\varepsilon}G(\lambda,\tau)d\tau\right)\overline{\Psi_{\alpha,\beta}^M(x,\lambda)}d\sigma(\varepsilon),$$

The proof is completed.
\section{Conclusion}

In this paper, we propose a novel harmonic analysis method, the Linear Canonical Jacobi-Dunkl Transform (LCJDT), by combining the strengths of the LCT and JDT. We begin by defining the linear canonical Jacobi-Dunkl operator, deriving its kernel function, and exploring several fundamental properties. The theory of the LCJDT is then thoroughly examined, including the derivation of its inversion formula, Parseval's identity, differentiation, convolution, and uncertainty principle. This analysis demonstrates that the LCJDT provides a powerful tool for addressing complex problems involving non-stationary signals, asymmetric structures, and weighted functions, outperforming classical integral transforms such as the FT and LT in these contexts.

Moreover, we explored the application of the LCJDT in solving the heat equation, highlighting its capability to address partial differential equations (PDE) with complex boundary conditions and asymmetric weight functions. The LCJDT proves to be particularly effective in scenarios where classical transforms such as the FT and LT fall short due to their limited flexibility in handling non-uniform boundary conditions and weighted domains. By introducing five degrees of freedom, the LCJDT not only enhances its adaptability but also offers precise control over the transformation process, making it well-suited for problems that require capturing intricate local features across different scales and orientations. This flexibility is particularly valuable in fields such as image encryption, where security algorithms demand robust directional feature extraction, and target detection, where accurate and multi-dimensional signal analysis is critical. The ability of the LCJDT to operate effectively across a wide range of applications underscores its potential as a versatile mathematical tool in solving advanced problems in applied mathematics, physics, and engineering.

several aspects still require further investigation to enhance its applicability and efficiency. One critical area for future research is the in-depth analysis of the splitting properties of the kernel function, which could provide a more granular understanding of how the LCJDT behaves under different boundary conditions and parameter configurations \cite{r19}. This investigation is particularly important for extending the transform's use in higher-dimensional and more intricate physical models, such as those encountered in fluid dynamics, quantum mechanics, and electromagnetics \cite{r50,r51}. Furthermore, the development of fast algorithms tailored to the LCJDT remains an open challenge. Efficient computational methods, such as those leveraging numerical techniques like the fast Fourier transform (FFT) or parallel processing strategies, could significantly improve the transform's performance in real-time applications, particularly in areas such as signal processing, image analysis, and machine learning \cite{r52}. Addressing these challenges will not only broaden the scope of the LCJDT's applications but also enhance its viability as a robust tool in modern computational mathematics \cite{r50}.

\textbf{Data availability}

No data was used for the research described in this paper.

\textbf{Acknowledgments}

This work was supported in part by the National Natural Science Foundation of China [No. 62171041]; the Natural Science Foundation of Beijing Municipality [No. 4242011].

\textbf{Declarations}

\textbf{Conflict of interest:} The authors declare that they have no known competing financial interests or personal relationships that could have appeared to influence the work reported in this paper.

\end{document}